# On Comparison of Two Reliable Techniques for the Riesz Fractional Complex Ginzburg-Landau-Schrödinger Equation in Modelling Superconductivity

Asim Patra

*National Institute of Technology*
*Department of Mathematics*
*Rourkela-769008, India*
*Email: asimp1993@gmail.com*

**Abstract**

In the present paper, the Complex Ginzburg-Landau-Schrödinger (CGLS) equation with the Riesz fractional derivative has been treated by a reliable implicit finite difference method (IFDM) of second order and furthermore for the purpose of a comparative study, and also for the investigation of the accuracy of the resulting solutions another effective spectral technique viz. time-splitting Fourier spectral (TSFS) technique has been utilized. In the case of the finite difference discretization, the Riesz fractional derivative is approximated by the fractional centered difference approach. Further the stability of the proposed methods has been analysed thoroughly and the TSFS technique is proved to be unconditionally stable. Moreover the absolute errors for the solutions of $|\psi(x,t)|^2$ obtained from both the techniques for various fractional order have been tabulated. Further the $L^2$ and $L^\infty$ error norms has been displayed for $|\psi(x,t)|^2$ and the results are also graphically depicted.

**Keywords:** Complex Ginzburg-Landau-Schrödinger (CGLS) equation, Riesz fractional derivative, fractional centered difference, time-splitting Fourier spectral method.

**PACS**: 02.30.Jr, **MSC code**: 26A33

1. **Introduction**

The concept of fractional differential equations is very old concept but in the recent years it had shown tremendous growth of interest among the researchers. The field of fractional calculus plays a significant role in modelling various anomalous transport phenomena. There are several applications of the fractional differential equations in many areas like engineering, physics, astrophysics, chemistry, control theory, dynamics of complex systems, fluid dynamics, other sciences and historical summaries for the growth and advancement of fractional calculus[1-4]. The concept of Riesz derivative can be explained as the mathematical combination of the the left and right Riemann-Lioville derivatives and it is utilized in many real life processes. Technically, the Riesz derivative has a dual and nonlocal nature for which it has become a very essential tool to be utilized in the fractional modeling viz. diffusion equation. The Riesz fractional derivative and its generalizations are used in several equations which gives a description regarding the applications in random walk models and anomalous diffusion characterized by nonlinear dependence of the mean square displacement of a diffusing particle over time.

The present paper deals with the following Riesz fractional Complex Ginzburg-Landau-Schrödinger(CGLS) equation

$$(\eta - i\beta)\frac{\partial \psi}{\partial t} = \frac{\partial^{\alpha} \psi}{\partial |x|^{\alpha}} + \frac{1}{\varepsilon^2}\left(V(x) - |\psi|^2\right)\psi, \quad x \in \Re, \; t > 0, \tag{1.1}$$

with initial condition

$$\psi(x,0) = \psi_0(x), \quad x \in \Re, \tag{1.2}$$

where $\psi(x,t)$ is a complex valued wave function, $\frac{\partial^{\alpha} \psi}{\partial |x|^{\alpha}}$ represents the Riesz derivative with fractional order $\alpha$ which is also denoted as $\frac{\partial^{\alpha} \psi}{\partial |x|^{\alpha}} = -(-\Delta)^{\alpha/2}\psi(x,t)$. Here $\eta$ and $\beta$ are two nonnegative constants satisfying $\eta + \beta > 0$, $\beta \neq 0$ in the complex case and the real-valued external potential $V(x)$ takes the form

$$V(x) = \frac{1}{1 + \exp(-\gamma_x x^2)},$$

where $\gamma_x$ is a positive constant.

Moreover $\varepsilon$ is a positive constant which takes the value unity in the standard CGLS equation but $0 < \varepsilon < 1$ for the CGLS equation in the semi-classical regime. Here the complex order parameter case is taken into consideration which defines the modelling of superconductivity process.

In the present article, the Riesz fractional Complex Ginzburg-Landau-Schrödinger (CGLS) equation has been treated by an effective scheme known as the time-splitting Fourier spectral method which can be otherwise described as a scheme of finding the solution in smaller steps. This technique is actually very efficient and powerful in handling the non-linear partial differential equations [5-7]. Apart from that an implicit finite difference method (IFDM) has been used for solving the Riesz fractional partial differential equation where the Riesz fractional derivative is approximated by the fractional centered difference scheme.

The literature regarding the numerical methodology of the Riesz fractional Complex Ginzburg-Landau-Schrödinger equation is limited to some case. Therefore, the techniques represented in the paper would prove to be highly desirable and effective. First of all the implicit finite difference discretization technique involved by Celik et al[8] for the Riesz fractional diffusion equation is included here. Lin et al[9] have depicted the dynamics, pinning and hyteresis of Ginzburg-Landau equation. Further Zhang et al[10] have studied the numerical simulation of vortex dynamics in the Ginzburg-Landau-Schrödinger equation. Jian et al[11] have presented the vortex dynamics in the Ginzburg-Landau equation in the inhomogeneous superconductors Bai and Wang[12] utilized the time splitting spectral scheme to the Schrödinger-Boussinesq equations and performed numerical experiments to show the numerical accuracy offered by the proposed method. Further, Markowich et al [13] have studied the numerical approximation of quadratic observables of Schrödinger type equations in the semi-classical limit. Wang [14] utilized the time splitting spectral method for the coupled G-P equations.

This organization of this research article is done in the following manner: firstly, some of the basic preliminaries regarding the Riesz fractional derivative and discrete Fourier transform have been presented in section 2. Then, in section 3, we describe the numerical implicit finite difference method for the Riesz fractional CGLS equation. Then stability analysis has been presented for the proposed implicit numerical method in section 4 via von Neumann stability analysis. In section 5, an effective fourier spectral approach is exhibited for the proposed fractional equation via splitting technique. Then in section 6, the splitting scheme is proved to be unconditionally stable. Furthermore in section 7, the numerical experiments and discussions have been presented and along with that, the tabulated results have been exhibited for the various solutions of absolute value of wave function. Furthermore, graphical solutions have been also exhibited here. Then section 8 is devoted for the conclusion.

**2. Basic preliminaries of fractional calculus**

In this section some of the basic concepts regarding the Riesz derivative have been discussed. Moreover some of the formulae, which will be later used, regarding the discrete Fourier transform and inverse discrete Fourier transform has been displayed.

**Definition 1**: For a fractional order $\alpha$ $(n-1<\alpha \leq n)$, the Riesz derivative on the infinite domain $-\infty < x < \infty$ is defined as [15]

$$\frac{\partial^\alpha \psi(x)}{\partial |x|^\alpha} = -c_\alpha \left( {}_{-\infty}D_x^\alpha \psi(x) + {}_xD_\infty^\alpha \psi(x) \right), \tag{2.1}$$

where ${}_{-\infty}D_x^\alpha \psi(x)$ is the Riemann-Liouville left derivative defined as

$$_{-\infty}D_x^\alpha \psi(x) = \frac{1}{\Gamma(n-\alpha)} \frac{\partial^n}{\partial x^n} \int_{-\infty}^x \frac{\psi(\xi)d\xi}{(x-\xi)^{1-n+\alpha}}, \tag{2.2}$$

and ${}_xD_\infty^\alpha \psi(x)$ is the Riemann-Liouville right derivative defined as

$$_xD_\infty^\alpha \psi(x) = \frac{(-1)^n}{\Gamma(n-\alpha)} \frac{\partial^n}{\partial x^n} \int_x^\infty \frac{\psi(\xi)d\xi}{(\xi-x)^{1-n+\alpha}}, \tag{2.3}$$

$$c_\alpha = \frac{1}{2\cos\frac{\alpha\pi}{2}}, \ \alpha \neq 1. \tag{2.4}$$

When $\alpha = 1$, we have

$$D_x^1 \psi(x) = \frac{d}{dx} H\psi(x) = \frac{d}{dx} \frac{1}{\pi} \int_{-\infty}^\infty \frac{\psi(z)dz}{z-x}, \tag{2.5}$$

where $H$ is the Hilbert transform and the integral is taken to be in the Cauchy principal value sense.

If $\psi(x,t)$ is given on the finite domain $[a,b]$ and it satisfies the boundary conditions $u(a,t)=u(b,t)=0$, the function can be extended by taking $\psi(x,t) \equiv 0$ for $x \leq a$ and $x \geq b$. In view of this, the Riesz fractional derivative of order $\alpha$ $(n-1<\alpha \leq n)$ on the finite interval $a \leq x \leq b$ can be defined as

$$\frac{\partial^{\alpha}\psi(x,t)}{\partial |x|^{\alpha}} = -\frac{1}{2\cos\frac{\alpha\pi}{2}}\left(_aD_x^{\alpha}\psi(x,t)+{}_xD_b^{\alpha}\psi(x,t)\right), \qquad (2.6)$$

where

$$_aD_x^{\alpha}\psi(x,t) = \frac{1}{\Gamma(n-\alpha)}\frac{\partial^n}{\partial x^n}\int_a^x \frac{\psi(\xi,t)d\xi}{(x-\xi)^{1-n+\alpha}}, \qquad (2.7)$$

$$_xD_b^{\alpha}\psi(x,t) = \frac{(-1)^n}{\Gamma(n-\alpha)}\frac{\partial^n}{\partial x^n}\int_x^b \frac{\psi(\xi,t)d\xi}{(\xi-x)^{1-n+\alpha}}. \qquad (2.8)$$

**Definition 3:** Now for solving the equation in the later section of the paper, some of the mathematical definitions should be known regarding the discrete Fourier transform applied here.

The discrete Fourier transformation of a sequence $\{\psi_j\}$ is defined from $t=t_n$ to $t_{n+1}$ as

$$\hat{\psi}_k(t) = F_k[\psi_j(t)] = \sum_{j=0}^{m-1}\psi_j(t)\exp(-i\mu_k(x_j-a)) \quad , \quad k=-\frac{m}{2},\ldots,\frac{m}{2}-1 \qquad (2.9)$$

where $\mu_k = \frac{2\pi k}{b-a}$ is the transform parameter and the domain of $x$ is $[a,b]$.

Similarly the formula for the inverse discrete Fourier transform is given by

$$\psi_j(t) = F_j^{-1}[\hat{\psi}_k(t)] = \frac{1}{m}\sum_{k=-\frac{m}{2}}^{\frac{m}{2}-1}\hat{\psi}_k(t)\exp(i\mu_k(x_j-a)) \quad , \quad j=0,1,2,\ldots,m-1. \qquad (2.10)$$

## 3. Second-order implicit finite difference method for the fractional Complex Ginzburg-Landau-Schrödinger equation with the Riesz fractional derivative

In this section an implicit finite difference technique has been utilized for dealing with the following Riesz fractional Complex Ginzburg-Landau-Schrödinger equation

$$(\eta-i\beta)\frac{\partial\psi}{\partial t} = -(-\Delta)^{\alpha/2}\psi(x,t) + \frac{1}{\varepsilon^2}\left(V(x)-|\psi|^2\right)\psi, \quad x\in\mathfrak{R},\ t>0, \qquad (3.1)$$

with initial condition

$$\psi(x,0) = \mu(x), \quad x\in\mathfrak{R}, \qquad (3.2)$$

where $\psi(x,t)$ is a complex valued wave function.

**Lemma 3.1**: Let $f\in C^5(\mathbb{R})$ with all the derivatives up to fifth order belong to the space $L_1(\mathbb{R})$ and the fractional centered difference be

$$\Delta_h^{\alpha}f(x) = \sum_{j=-\infty}^{\infty}\frac{(-1)^j\Gamma(\alpha+1)}{\Gamma\left(\frac{\alpha}{2}-j+1\right)\Gamma\left(\frac{\alpha}{2}+j+1\right)}f(x-jh), \qquad (3.3)$$

then

$$-h^{-\alpha}\Delta_h^{\alpha}f(x) = \frac{\partial^{\alpha}f(x)}{\partial |x|^{\alpha}} + O(h^2). \qquad (3.4)$$

So here when $h \to 0$, $\dfrac{\partial^\alpha f(x)}{\partial |x|^\alpha}$ represents the Riesz derivative of fractional order $\alpha$ for $1 < \alpha \leq 2$.

Proof: It may be referred to Ref. [8] for proof of the above Lemma.

**Corollary**:

Recently, in Ref. [8], Çelik and Duman has derived the interesting result if $f^*(x)$ is defined by

$$f^*(x) = \begin{cases} f(x), & x \in [a,b] \\ 0, & x \notin [a,b] \end{cases}$$

such that $f^* \in C^5(\mathbb{R})$ and all derivatives upto fifth order belongs to the space $L_1(\mathbb{R})$, then for the Riesz derivative of fractional order $\alpha$ $(1 < \alpha \leq 2)$

$$\frac{\partial^\alpha f(x)}{\partial |x|^\alpha} = -h^{-\alpha} \sum_{j=-\frac{b-x}{h}}^{\frac{x-a}{h}} \frac{(-1)^j \Gamma(\alpha+1)}{\Gamma\left(\frac{\alpha}{2} - j + 1\right)\Gamma\left(\frac{\alpha}{2} + j + 1\right)} f(x - jh) + O(h^2), \quad (3.5)$$

where $h = \dfrac{b-a}{m}$, and $m$ denotes the number of subintervals of the interval $[a,b]$.

### 3. 1. Second-order implicit finite difference approach for the CGLS equation involving Riesz fractional derivative by the fractional centered difference

Now the second-order implicit finite difference discretization for the equation (3.1) is shown as

$$(\eta - i\beta)\left(\frac{\psi_i^{n+1} - \psi_i^n}{\tau}\right) - \frac{1}{2}(-h^{-\alpha})\left(\sum_{j=i-m}^{i} g_j \psi_{i-j}^{n+1} + \sum_{j=i-m}^{i} g_j \psi_{i-j}^n\right)$$

$$= \frac{1}{\varepsilon^2}\left(V(x)\left(\frac{\psi_i^{n+1} + \psi_i^n}{2}\right) - \frac{1}{2}\left(\left(\psi_i^{n+1}\right)^2 \overline{\psi}_i^{n+1} + \left(\psi_i^n\right)^2 \overline{\psi}_i^n\right)\right) + R_i^{n+\frac{1}{2}}, \quad (3.6)$$

where $R_i^{n+\frac{1}{2}}$ is the local truncation error given by $R_i^{n+\frac{1}{2}} = O(\tau^2 + h^2)$.

### 4. Von Neumann stability analysis for the proposed implicit finite difference scheme

**Theorem 4.1:** The proposed implicit finite difference scheme in eq. (5.4) is unconditionally stable.

**Proof:** The proposed implicit scheme is given by

$$\psi_i^{n+1} = \psi_i^n - \frac{\tau}{2(\eta - \sqrt{-1}\beta)} h^{-\alpha} \sum_{j=i-m}^{i} g_j \psi_{i-j}^{n+1} - \frac{\tau}{2(\eta - \sqrt{-1}\beta)} h^{-\alpha} \sum_{j=i-m}^{i} g_j \psi_{i-j}^n$$

$$+ \frac{\tau}{(\eta - \sqrt{-1}\beta)\varepsilon^2}\left(V(x)\left(\frac{\psi_i^{n+1} + \psi_i^n}{2}\right) - \frac{1}{2}\left(\left(\psi_i^{n+1}\right)^2 \overline{\psi}_i^{n+1} + \left(\psi_i^n\right)^2 \overline{\psi}_i^n\right)\right) \quad (4.1)$$

where $\beta \neq 0$.

Let $\psi_i^n = \xi^n e^{\sqrt{-1}\omega i h}$, with $\xi > 0$. Then, we obtain

$$\xi^{n+1} e^{\sqrt{-1}\omega i h} = \xi^n e^{\sqrt{-1}\omega i h} - \frac{\tau}{h^\alpha (\eta - \sqrt{-1}\beta)}\left(\frac{1}{2}\sum_{j=i-m}^{i} g_j \xi^{n+1} e^{\sqrt{-1}\omega(i-j)h} + \frac{1}{2}\sum_{j=i-m}^{i} g_j \xi^n e^{\sqrt{-1}\omega(i-j)h}\right)$$

$$+\frac{\tau}{2(\eta-\sqrt{-1}\beta)\varepsilon^2}\left(V(x)\left(\xi^{n+1}e^{\sqrt{-1}\omega ih}+\xi^{n}e^{\sqrt{-1}\omega ih}\right)-(\xi^{n+1}e^{\sqrt{-1}\omega ih})^2\xi^{n+1}e^{-\sqrt{-1}\omega ih}-(\xi^{n}e^{\sqrt{-1}\omega ih})^2\xi^{n}e^{-\sqrt{-1}\omega ih}\right) \quad (4.2)$$

This can also be written as

$$\xi^{n+1}e^{\sqrt{-1}\omega ih}=\xi^{n}e^{\sqrt{-1}\omega ih}-\frac{r}{(\eta-\sqrt{-1}\beta)}\left(\sum_{j=i-m}^{i}g_j\xi^{n+1}e^{\sqrt{-1}\omega(i-j)h}+\sum_{j=i-m}^{i}g_j\xi^{n}e^{\sqrt{-1}\omega(i-j)h}\right)$$

$$+\frac{\tau}{2(\eta-\sqrt{-1}\beta)\varepsilon^2}\left(V(x)\left(\xi^{n+1}e^{\sqrt{-1}\omega ih}+\xi^{n}e^{\sqrt{-1}\omega ih}\right)-(\xi^{n+1}e^{\sqrt{-1}\omega ih})^2\xi^{n+1}e^{-\sqrt{-1}\omega ih}-(\xi^{n}e^{\sqrt{-1}\omega ih})^2\xi^{n}e^{-\sqrt{-1}\omega ih}\right), \quad (4.3)$$

where $r=\dfrac{\tau}{2h^\alpha}$.

Let us assume that

$$\Psi_{\max}^{(n)}=\max_{0\leq i\leq m-1}\left\{\left(\Psi_i^n\right)^2,\left(\Psi_i^{n+1}\right)^2\right\}, \quad (4.4)$$

$$\xi^{n+1}e^{\sqrt{-1}\omega ih}\left(1+\frac{r}{(\eta-\sqrt{-1}\beta)}\sum_{j=i-m}^{i}g_je^{-\sqrt{-1}\omega jh}-\frac{\tau}{2(\eta-\sqrt{-1}\beta)}V(x)\right)+\frac{\tau}{2(\eta-\sqrt{-1}\beta)\varepsilon^2}\Psi_{\max}^{(n)}\xi^{n+1}e^{\sqrt{-1}\omega ih}$$

$$=\xi^{n}e^{\sqrt{-1}\omega ih}\left(1-\frac{r}{(\eta-\sqrt{-1}\beta)}\sum_{j=i-m}^{i}g_je^{-\sqrt{-1}\omega jh}+\frac{\tau}{2(\eta-\sqrt{-1}\beta)\varepsilon^2}V(x)\right)-\frac{\tau}{2(\eta-\sqrt{-1}\beta)\varepsilon^2}\Psi_{\max}^{(n)}\xi^{n}e^{-\sqrt{-1}\omega ih} \quad (4.5)$$

Now further taking the absolute of both the sides and triangle inequality,

$$\left|\xi^{n+1}e^{\sqrt{-1}\omega ih}\left(1+\frac{r}{(\eta-\sqrt{-1}\beta)}\sum_{j=i-m}^{i}g_je^{-\sqrt{-1}\omega jh}-\frac{\tau}{2(\eta-\sqrt{-1}\beta)}V(x)\right)\right|$$

$$\leq\left|\xi^{n}e^{\sqrt{-1}\omega ih}\left(1-\frac{r}{(\eta-\sqrt{-1}\beta)}\sum_{j=i-m}^{i}g_je^{-\sqrt{-1}\omega jh}+\frac{\tau}{2(\eta-\sqrt{-1}\beta)\varepsilon^2}V(x)\right)-\frac{\tau}{2(\eta-\sqrt{-1}\beta)\varepsilon^2}\Psi_{\max}^{(n)}\xi^{n}e^{-\sqrt{-1}\omega ih}\right|+\left|\frac{\tau}{2(\eta-\sqrt{-1}\beta)\varepsilon^2}\Psi_{\max}^{(n)}\xi^{n+1}e^{\sqrt{-1}\omega ih}\right|$$

(4.6)

After eliminating $\xi^n e^{\sqrt{-1}\omega ih}$ from the eq. (4.6) and rearranging the terms, we get

$$|\xi|\left(\left|1+\frac{r}{(\eta-\sqrt{-1}\beta)}\sum_{j=i-m}^{i}g_je^{-\sqrt{-1}\omega jh}-\frac{\tau}{2(\eta-\sqrt{-1}\beta)}V(x)\right|-\left|\frac{\tau}{2(\eta-\sqrt{-1}\beta)\varepsilon^2}\Psi_{\max}^{(n)}\right|\right)$$

$$\leq\left|\left(1-\frac{r}{(\eta-\sqrt{-1}\beta)}\sum_{j=i-m}^{i}g_je^{-\sqrt{-1}\omega jh}+\frac{\tau}{2(\eta-\sqrt{-1}\beta)\varepsilon^2}V(x)\right)-\frac{\tau}{2(\eta-\sqrt{-1}\beta)\varepsilon^2}\Psi_{\max}^{(n)}\right| \quad (4.7)$$

Now by the triangle inequality we know that,

$$|a_1|-|a_2|\leq|a_1-a_2|,$$

This implies that

$$\left|1+\frac{r}{(\eta-\sqrt{-1}\beta)}\sum_{j=i-m}^{i}g_j e^{-\sqrt{-1}\omega jh}-\frac{\tau}{2(\eta-\sqrt{-1}\beta)}V(x)\right|-\left|\frac{\tau}{2(\eta-\sqrt{-1}\beta)\varepsilon^2}\Psi_{\max}^{(n)}\right|$$

$$\leq\left|1+\frac{r}{(\eta-\sqrt{-1}\beta)}\sum_{j=i-m}^{i}g_j e^{-\sqrt{-1}\omega jh}-\frac{\tau}{2(\eta-\sqrt{-1}\beta)}V(x)-\frac{\tau}{2(\eta-\sqrt{-1}\beta)\varepsilon^2}\Psi_{\max}^{(n)}\right|. \qquad (4.8)$$

So by combining eq. (4.7) and eq. (4.8), we observe that in order to satisfy both the inequalities we require

$$|\xi|\leq 1.$$

So, it has been proved that the proposed implicit finite difference scheme is unconditionally stable for the fractional CGLS equation as per the Von Neumann's stability analysis.  ∎

## 5. The proposed spectral method for the fractional Ginzburg-Landau-Schrödinger equation with the Riesz fractional derivative

An effective numerical scheme has been developed in this section for dealing with the following Riesz fractional Ginzburg-Landau-Schrödinger equation

$$(\eta-i\beta)\frac{\partial\psi}{\partial t}=-(-\Delta)^{\alpha/2}\psi(x,t)+\frac{1}{\varepsilon^2}\left(V(x)-|\psi|^2\right)\psi, \quad x\in\Re, \ t>0, \qquad (5.1)$$

with initial condition

$$\psi(x,0)=\mu(x), \quad x\in\Re, \qquad (5.2)$$

where $\psi(x,t)$ is a complex valued wave function.

### i. Discretisation procedure

The mesh size is $h=\frac{b-a}{m}$, $m$ denoting an even integer and $x_j=a+jh$ is taken as the grid points where $j=0,1,...,m$. Similarly taking temporal step size $\tau$, the time steps taken as $t_n=n\tau$, $\tau>0$, $n=0,1,2,3,...$

### ii. Methodology

The first Schrödinger–like equation is tackled in two steps, rather splitting from $t=t_n$ to $t_{n+1}$. So it is split in the following way by first solving

$$(\eta-i\beta)\frac{\partial\psi}{\partial t}=-(-\Delta)^{\alpha/2}\psi(x,t), \qquad (5.3)$$

and then solving

$$(\eta-i\beta)\frac{\partial\psi}{\partial t}=\frac{1}{\varepsilon^2}\left(V(x)-|\psi|^2\right)\psi, \qquad (5.4)$$

by taking time step $\tau$ in both the equations.

So by applying the spectral technique via TSFS approach for discretizing eq. (5.3) in spatial variable, along with applying an integration for eq. (5.4) in temporal variable by taking integrating limits from $t_n$ to $t_{n+1}$ we obtain

$$\psi(x, t_{n+1}) = \exp\left[\int_{t_n}^{t_{n+1}} \left(\frac{1}{\varepsilon^2(\eta - i\beta)}(V(x) - \psi(x,s)\overline{\psi}(x,s))\right)ds\right]\psi(x, t_n) \qquad (5.5)$$

by applying a numerical method via trapezoidal formula for the purpose of approximation of the integral on $[t_n, t_{n+1}]$ which provides us with the solution

$$\psi(x, t_{n+1}) = \exp\left[\frac{1}{2\varepsilon^2(\eta - i\beta)}\tau(2V(x) - (\psi(x, t_{n+1})\overline{\psi}(x, t_{n+1}) + \psi(x, t_n)\overline{\psi}(x, t_n))\right]\psi(x, t_n). \qquad (5.6)$$

Now solving equation (5.3) by applying the Fourier transform technique to

$$(\eta - i\beta)\frac{\partial \psi}{\partial t} = -(-\Delta)^{\alpha/2}\psi(x, t), \qquad (5.7)$$

and further integrating it from $t = t_n$ to $t_{n+1}$, we obtain

$$\hat{\psi}(\mu_k, t_{n+1}) = \exp\left(-\frac{1}{\eta - i\beta}|\mu_k|^\alpha \tau\right)\hat{\psi}(\mu_k, t_n), \qquad (5.8)$$

where $\mu_k$ is the transform parameter taken as $\mu_k = \frac{2\pi k}{b-a}$.

Taking the inverse discrete Fourier transform to eq. (5.8), we obtain the simplified solution as

$$\psi_j^{**} = \frac{1}{m}\sum_{k=-\frac{m}{2}}^{\frac{m}{2}-1} \exp(-\frac{1}{\eta - i\beta}|\mu_k|^\alpha \tau)(\hat{\psi}^*)_k \exp(i\mu_k(x_j - a)). \qquad (5.9)$$

where $(\hat{\psi}^*)_k$ is the Fourier transform of $\psi_j^*$ defined by

$$(\hat{\psi}^*)_k(t) = F_k[\psi_j^*(t)] = \sum_{j=0}^{m-1}\psi_j^*(t)\exp(-i\mu_k(x_j - a)) \quad , \quad k = -\frac{m}{2}, ..., \frac{m}{2} - 1. \qquad (5.10)$$

### iii. Basic idea of the methodology of the second order splitting scheme via Strang splitting technique

If an ordinary differential equation can be written as

$$f' = L_1(f) + L_2(f),$$

where the exact solutions of $f' = L_1(f)$, denoted by $\phi_t^{L_1}(f(0))$ and $f' = L_2(f)$, denoted by $\phi_t^{L_2}(f(0))$ can be efficiently computed, the splitting techniques often provide an effective alternative in comparison with the traditional integration methods. Moreover the splitting schemes preserve a given property of the ordinary differential equation which is done by the flow generated by $L_1$ and $L_2$. Now if the exact partial flows are utilized the for a step size $\tau$, the Strang splitting scheme can be manifested as

$$S_\tau = \phi^{L_1}_{\frac{\tau}{2}} \circ \phi^{L_2}_\tau \circ \phi^{L_1}_{\frac{\tau}{2}},$$

which is a symmetric technique of second order.

The second order scheme can be otherwise implemented in the following steps

$$\tilde{f}_1 = e^{L_1\frac{\tau}{2}} y_0, \quad \vec{f}_1 = e^{L_2\tau} \tilde{f}_1, \quad f_1 = e^{L_1\frac{\tau}{2}} \vec{f}_1,$$

$$\tilde{f}_2 = e^{L_1\frac{\tau}{2}} f_1, \quad \vec{f}_2 = e^{L_2\tau} \tilde{f}_2, \quad f_2 = e^{L_1\frac{\tau}{2}} \vec{f}_2,$$

$$\tilde{f}_3 = e^{L_1\frac{\tau}{2}} f_2, \quad \vec{f}_3 = e^{L_2\tau} \tilde{f}_3, \quad f_3 = e^{L_1\frac{\tau}{2}} \vec{f}_3,$$

.

.

$$\tilde{f}_n = e^{L_1\frac{\tau}{2}} f_{n-1}, \quad \vec{f}_n = e^{L_2\tau} \tilde{f}_n, \quad y_n = e^{L_1\frac{\tau}{2}} \vec{f}_n.$$

**iv. Implementation of the Strang splitting (SP) scheme for the proposed Riesz fractional GLS equation**

The technique of Strang splitting scheme is displayed by the following three equations

$$\psi^*_j = \exp\left(\frac{(\eta+i\beta)}{2(\eta^2+\beta^2)\varepsilon^2}\tau(2V(x)-(\psi(x,t_{n+1})\overline{\psi}(x,t_{n+1})+\psi(x,t_n)\overline{\psi}(x,t_n)))\right)\psi(x,t_n) \quad, \quad j=0,1,2,...,m-1,$$

(5.11)

$$\psi^{**}_j = \frac{1}{m}\sum_{k=-\frac{m}{2}}^{\frac{m}{2}-1} \exp(-\frac{(\eta+i\beta)}{(\eta^2+\beta^2)}|\mu_k|^\alpha \tau)(\hat{\psi}^*)_k \exp(i\mu_k(x_j-a)) \quad, \quad j=0,1,2,...,m-1, \qquad (5.12)$$

$$\psi^{n+1}_j = \exp\left(\frac{(\eta+i\beta)}{2(\eta^2+\beta^2)\varepsilon^2}\tau(2V(x)-(\psi(x,t_{n+1})\overline{\psi}(x,t_{n+1})+\psi(x,t_n)\overline{\psi}(x,t_n)))\right)\psi^{**}_j \quad, \quad j=0,1,2,...,m-1,$$

(5.13)

where $(\hat{\psi}^*)_k$ is the of Fourier transform $\psi^*_j$ as the earlier given definition.

**6. Stability analysis for the proposed spectral scheme for the CGLS equation involving Riesz fractional derivative**

**Theorem 6.1**  The second-order time-splitting scheme (5.11), (5.12) and (5.13) for the Riesz fractional Ginzburg-Landau-Schrödinger equation is unconditionally stable.

Furthermore, the conservative properties are shown as

$$\|\psi^{n+1}\|^2_{l_2} = \|\psi^0\|^2_{l_2} \quad, \quad n=0,1,2,... \qquad (6.1)$$

**Proof.** First of all here the definition for the $L^2-norm$ and discrete $l^2-norm$ are defined as

$$\|\psi\|_{L^2} = \sqrt{\int_a^b |\psi(x)|^2 dx}, \quad \|\psi\|_{l^2} = \sqrt{\frac{b-a}{m}\sum_{j=0}^{m-1}|\psi_j|^2}. \tag{6.2}$$

For the scheme (5.11)-(5.13), using (5.9), (5.10) and (6.2), we have

$$\frac{1}{b-a}\|\psi^{n+1}\|_{l_2}^2 = \frac{1}{m}\sum_{j=0}^{m-1}|\psi_j^{n+1}|^2$$

$$= \frac{1}{m}\sum_{j=0}^{m-1}\left|\exp\left(\frac{(\eta+i\beta)}{2(\eta^2+\beta^2)\varepsilon^2}\tau(2V(x)-(\psi(x,t_{n+1})\overline{\psi}(x,t_{n+1})+\psi(x,t_n)\overline{\psi}(x,t_n)))\right)\psi_j^{**}\right|^2 = \frac{1}{m}\sum_{j=0}^{m-1}|\psi_j^{**}|^2$$

$$= \frac{1}{m}\sum_{j=0}^{m-1}\left|\frac{1}{m}\sum_{k=-\frac{m}{2}}^{\frac{m}{2}-1}\exp(-\frac{(\eta+i\beta)}{(\eta^2+\beta^2)}|\mu_k|^\alpha\tau)(\hat{\psi}^*)_k \exp(i\mu_k(x_j-a))\right|^2 = \frac{1}{m^2}\sum_{k=-\frac{m}{2}}^{\frac{m}{2}-1}\left|\exp\left(-\frac{(\eta+i\beta)}{(\eta^2+\beta^2)}|\mu_k|^\alpha\tau\right)(\hat{\psi}^*)_k\right|^2$$

$$= \frac{1}{m^2}\sum_{k=-\frac{m}{2}}^{\frac{m}{2}-1}|(\hat{\psi}^*)_k|^2 = \frac{1}{m^2}\sum_{k=-\frac{m}{2}}^{\frac{m}{2}-1}\left|\sum_{j=0}^{m-1}\psi_j^* \exp(-i\mu_k(x_j-a))\right|^2 = \frac{1}{m}\sum_{j=0}^{m-1}|\psi_j^*|^2$$

$$= \frac{1}{m}\sum_{j=0}^{m-1}|\psi_j^n|^2 = \frac{1}{b-a}\|\psi^n\|_{l_2}^2. \tag{6.3}$$

Here the following identities have been used

$$\sum_{j=0}^{m-1}\exp(i2\pi(l-k)j/M) = \begin{cases}0, & l-k \neq mM \\ M, & l-k = mM\end{cases}, \quad M \text{ is an integer} \tag{6.4}$$

and

$$\sum_{k=-\frac{m}{2}}^{\frac{m}{2}-1}\exp(i2\pi(l-j)k/M) = \begin{cases}0, & l-j \neq mM \\ M, & l-j = mM\end{cases}, \quad M \text{ is an integer}. \tag{6.5}$$

Therefore, the result in eq. (6.1) can be obtained from (6.3) for the splitting scheme (5.11)-(5.13) by means of induction.

## 7. Numerical discussion

**Example 1** We consider the CGLS equation (3.1) with the initial condition

$$\psi_0(x) = \sqrt{1-\left(\frac{20\pi}{L}\right)^2}\exp\left(i\frac{20\pi}{L}x\right). \tag{7.1}$$

Here $L=100$ and the constants in eq. (3.1) are taken as $\eta =1$ and $\beta =1$ respectively. Here the interval for $x$ is taken as $[-5,5]$. Also the value for $\varepsilon =1$ and the constant in external potential $V(x)$ is taken as $\gamma_x =1$. So $V(x)= \dfrac{1}{1+\exp(-x^2)}$.

The section includes the numerical discussion and representation of the results obtained by the time-splitting Fourier spectral method for the CGLS equation (5.1) and (5.2). In this context, Table 1 contains the comparison for the values of $|\psi(x,t)|^2$, $\text{Re}(\psi(x,t))$ and $\text{Im}(\psi(x,t))$ by the TSFS scheme and IFDM scheme by providing the absolute errors for Example 1 for various values of $x$ for $-5\le x\le 5$ at $t=0.5$ for spatial step size $h=0.2$ and temporal step size $\tau =0.1$ taking two cases $\alpha =1.5$ and $\alpha =1.75$. Apart from that the error norms $L^2$–norm and $L^\infty$–norm for the absolute values, real values and imaginary values for the wave function have been displayed in Table 2 for various values of $x$ and $t$ for $\alpha =1.5$ and $\alpha =1.75$.

**Table 1** Absolute errors for fractional Complex Ginzburg-Landau-Schrödinger (CGLS) equation (5.1) obtained by using time-splitting Fourier spectral approximation and the implicit finite difference solutions for Example 1 at various points of $x$ and $t=0.5$ taking step size, $\tau =0.1$ and $\alpha =1.5$ and $\alpha =1.75$.

| | $\alpha=1.5$ (Example 1) | $\alpha=1.75$ (Example 1) | $\alpha=1.5$ (Example 1) | $\alpha=1.75$ (Example 1) | $\alpha=1.5$ (Example 1) | $\alpha=1.75$ (Example 1) |
|---|---|---|---|---|---|---|
| $x$ | $t=0.5$ | $t=0.5$ | $t=0.5$ | $t=0.5$ | $t=0.5$ | $t=0.5$ |
| | Absolute errors between $\|\psi\|^2{}_{TSFS}$ and $\|\psi\|^2{}_{IFDM}$ | Absolute errors between $\|\psi\|^2{}_{TSFS}$ and $\|\psi\|^2{}_{IFDM}$ | Absolute errors between $\text{Re}(\psi)_{TSFS}$ and $\text{Re}(\psi)_{IFDM}$ | Absolute errors between $\text{Re}(\psi)_{TSFS}$ and $\text{Re}(\psi)_{IFDM}$ | Absolute errors between $\text{Im}(\psi)_{TSFS}$ and $\text{Im}(\psi)_{IFDM}$ | Absolute errors between $\text{Im}(\psi)_{TSFS}$ and $\text{Im}(\psi)_{IFDM}$ |
| -4.8 | 1.0033E-4 | 1.2083E-4 | 3.2087E-3 | 1.4736E-2 | 5.6849E-2 | 8.1894E-2 |
| -4.0 | 1.0056E-5 | 2.5915E-7 | 1.4175E-2 | 8.1883E-2 | 2.16643E-2 | 3.0983E-3 |
| -3.4 | 3.3765E-7 | 3.7743E-5 | 3.1259E-2 | 9.6138E-2 | 5.9323E-3 | 5.0476E-2 |
| -2.8 | 3.4743E-5 | 3.9741E-5 | 2.5270E-2 | 2.4288E-3 | 4.7090E-2 | 6.7967E-3 |
| -2.0 | 1.3558E-4 | 2.6029E-4 | 2.8179E-3 | 1.0852E-2 | 6.2378E-2 | 7.0617E-2 |
| -1.2 | 3.4166E-4 | 2.5277E-4 | 7.6303E-3 | 2.6167E-2 | 6.1600E-2 | 3.5265E-2 |
| -0.4 | 5.7204E-5 | 5.4308E-4 | 2.9220E-2 | 8.72495E-3 | 6.15342E-3 | 1.7022E-2 |
| 0.4 | 2.68671E-4 | 5.6204E-5 | 7.2001E-3 | 9.7932E-3 | 3.4291E-2 | 4.3205E-3 |
| 1.2 | 1.5544E-4 | 1.3548E-4 | 1.1562E-2 | 1.5159E-2 | 3.1303E-2 | 3.5453E-3 |

| 2.0 | 3.30865E-4 | 2.018E-4 | 5.3281E-3 | 2.1847E-3 | 5.5893E-2 | 7.5154E-2 |
| 2.8 | 5.1080E-5 | 9.2081E-6 | 7.2196E-4 | 2.3859E-2 | 4.5733E-2 | 7.0903E-3 |
| 3.4 | 1.53906E-4 | 1.45906E-4 | 2.2517E-3 | 7.5222E-4 | 4.4443E-3 | 6.2556E-3 |
| 4.0 | 3.20545E-6 | 1.8843E-7 | 1.0017E-3 | 1.9566E-3 | 8.1083E-3 | 5.9469E-2 |
| 4.8 | 2.65729E-6 | 2.6095E-7 | 9.443E-3 | 7.028E-3 | 1.70455E-2 | 7.6215E-3 |

Now the error norms $L_2$ and $L_\infty$ are tabulated in Table 1. Here for a fixed time $t_l$, the $L_2$ norm for error is defined as

$$L_2 \equiv \sqrt{\frac{1}{2m+1} \sum_{i=1}^{2m+1} (\psi_{IFDM}(x_i,t_l) - \psi_{TSFS}(x_i,t_l))^2}$$

and $L_\infty$ norm for error is given by

$$L_\infty \equiv \underset{-m \leq i \leq m}{Max} |\psi_{IFDM}(x_i,t_l) - \psi_{TSFS}(x_i,t_l)|.$$

**Table 2** $L_2$ norm and $L_\infty$ norm of errors between the solutions of TSFS method and implicit finite difference method for Example 1 at $t = 0.5$ and 1, for $\alpha = 1.5$ and $\alpha = 1.75$.

| $t$ | Example 1 $\alpha = 1.5$ | | Example 1 $\alpha = 1.5$ | | Example 1 $\alpha = 1.5$ | |
|---|---|---|---|---|---|---|
| | $L_2$ error of $|\psi(x,t)|^2$ | $L_\infty$ error of $|\psi(x,t)|^2$ | $L_2$ error of Re($\psi(x,t)$) | $L_\infty$ error of Re($\psi(x,t)$) | $L_2$ error of Im($\psi(x,t)$) | $L_\infty$ error of Im($\psi(x,t)$) |
| 0.5 | 3.8237E-4 | 2.9388E-2 | 1.8793E-2 | 4.166E-2 | 4.735E-2 | 8.3122E-2 |
| 1.0 | 1.8071E-2 | 3.2506E-2 | 3.5968E-2 | 2.1704E-1 | 9.7688E-2 | 7.072E-2 |
| $t$ | Example 1 $\alpha = 1.75$ | | Example 1 $\alpha = 1.75$ | | Example 1 $\alpha = 1.75$ | |
| | $L_2$ error of $|\psi(x,t)|^2$ | $L_\infty$ error of $|\psi(x,t)|^2$ | $L_2$ error of Re($\psi(x,t)$) | $L_\infty$ error of Re($\psi(x,t)$) | $L_2$ error of Im($\psi(x,t)$) | $L_\infty$ error of Im($\psi(x,t)$) |
| 0.5 | 9.36681E-4 | 5.555E-3 | 1.80006E-2 | 4.0034E-2 | 4.794E-2 | 7.272E-4 |
| 1.0 | 2.7596E-3 | 3.9618E-3 | 2.149E-2 | 8.501E-2 | 3.5003E-2 | 7.7483E-2 |

Apart from the tabulated results, the solutions obtained by the proposed techniques have been graphically manifested by depicting the 2D as well as 3D plots for the solutions of $|\psi(x,t)|^2$,

Re($\psi(x,t)$)) and Im($\psi(x,t)$)) by taking fractional order $\alpha = 1.5$ and $\alpha = 1.75$ separately, for various values of $x$ in the interval $[-5, 5]$ and at $t = 0.5$. Apart from that comparison graphs obtained by both the TSFS and IFD techniques have also been depicted here for $|\psi(x,t)|^2$, Re($\psi(x,t)$) and Im($\psi(x,t)$) by taking fractional order $t = 0.5$ and $\alpha = 1.75$ which manifests fine tuning between the results. In Fig. 1, Fig. 2 and Fig. 3 the graphical 2D solutions for $|\psi(x,t)|^2$, Re($\psi(x,t)$) and Im($\psi(x,t)$) are depicted for Example 1 for both cases i.e. taking $\alpha = 1.5$ and $\alpha = 1.75$ separately. Further the corresponding 3D figures obtained by the time-splitting Fourier spectral method for $|\psi(x,t)|^2$, Re($\psi(x,t)$) and Im($\psi(x,t)$) for Example 1 are shown in Fig. 4, Fig. 5 and Fig. 6 respectively. Moreover the graphs obtained by the comparison of the TSFS technique and the IFDM are displayed in Fig. 7 for $|\psi(x,t)|^2$, in Fig. 8 for Re($\psi(x,t)$) and Im($\psi(x,t)$) in Fig. 9 by taking fractional order for two cases, $\alpha = 1.5$ and $\alpha = 1.75$. separately.

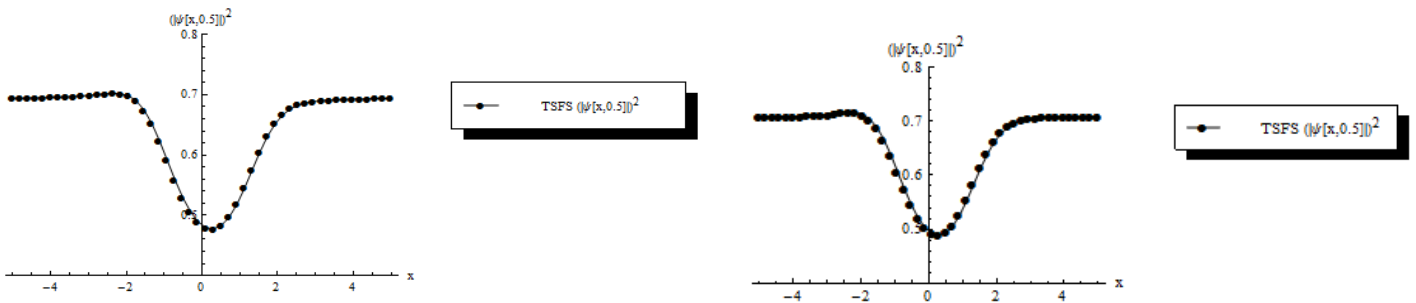

**Fig. 1.** Graphical TSFS 2D solutions for $|\psi(x,t)|^2$ of Example 1 for the Riesz fractional Complex Ginzburg-Landau-Schrödinger (CGLS) equation for various points of $x$ at $t = 0.5$ for $\alpha = 1.5$ and $\alpha = 1.75$.

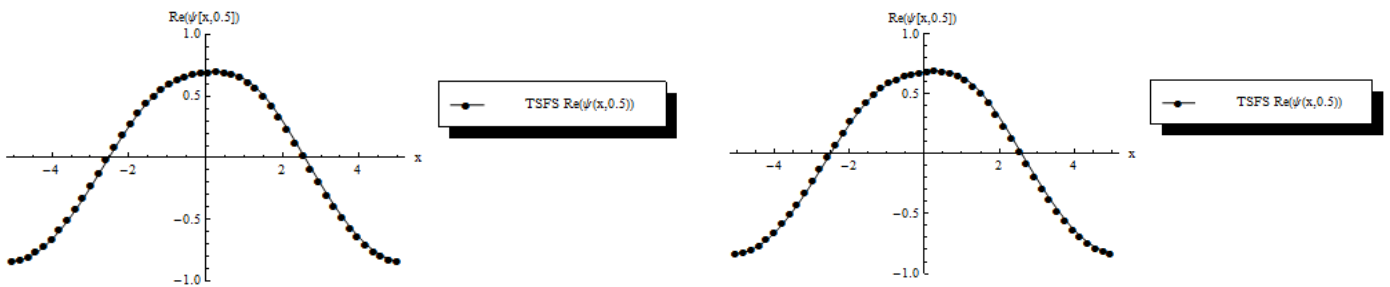

**Fig. 2.** Graphical TSFS 2D solutions for Re($\psi(x,t)$) of Example 1 for the Riesz fractional Complex Ginzburg-Landau-Schrödinger (CGLS) equation for various points of $x$ at $t = 0.5$ for $\alpha = 1.5$ and $\alpha = 1.75$.

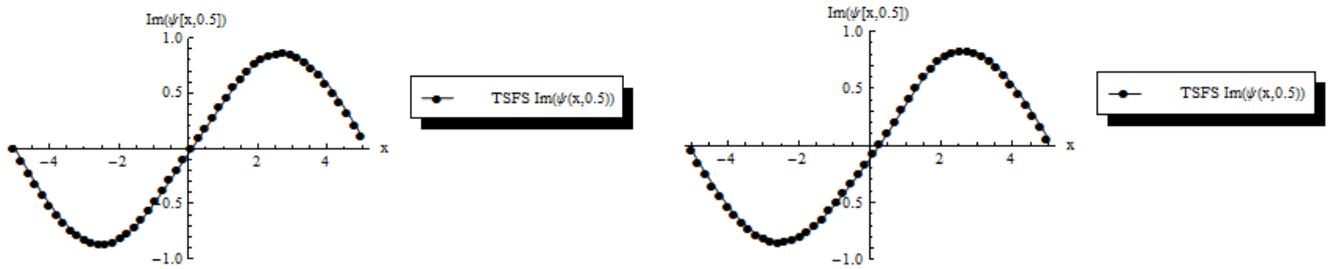

**Fig. 3.** Graphical TSFS 2D solutions for Im($\psi(x,t)$) of Example 1 for the Riesz fractional Complex Ginzburg-Landau-Schrödinger (CGLS) equation for various points of $x$ at $t = 0.5$ for $\alpha = 1.5$ and $\alpha = 1.75$.

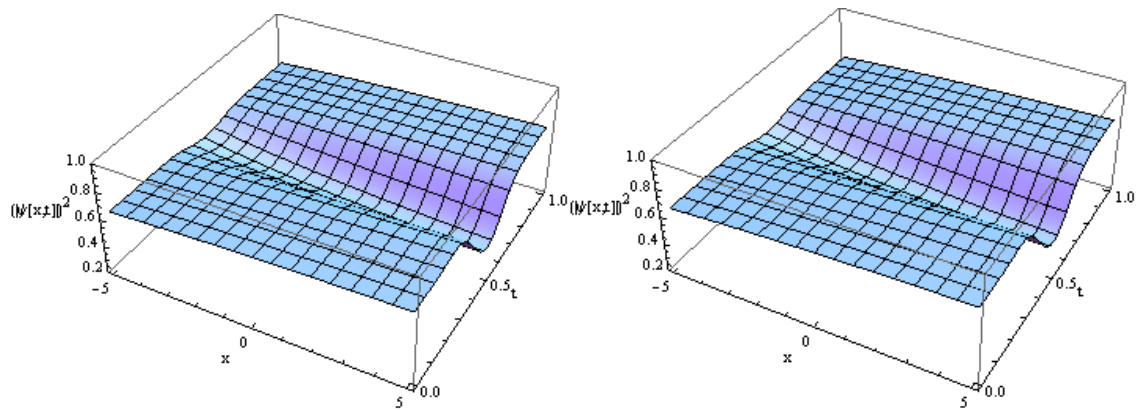

**Fig. 4.** Graphical TSFS 3D solutions for $|\psi(x,t)|^2$ of Example 1 for the Riesz fractional Complex Ginzburg-Landau-Schrödinger (CGLS) equation for various points of $x$ at $t = 0.5$ for $\alpha = 1.5$ and $\alpha = 1.75$.

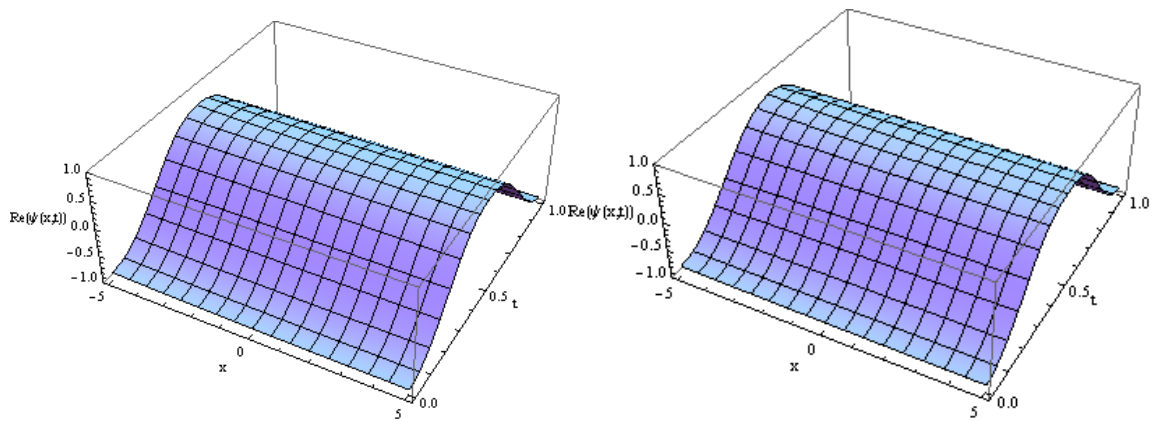

**Fig. 5.** Graphical TSFS 3D solutions for Re($\psi(x,t)$) of Example 1 for the Riesz fractional Complex Ginzburg-Landau-Schrödinger (CGLS) equation for various points of $x$ at $t = 0.5$ for $\alpha = 1.5$ and $\alpha = 1.75$.

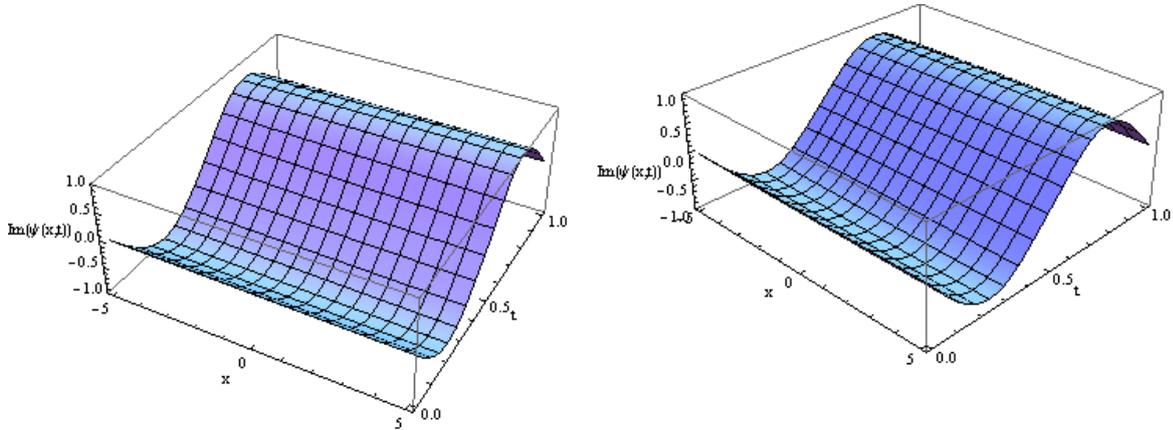

**Fig. 6.** Graphical TSFS 3D solutions for Im($\psi(x,t)$) of Example 1 for the Riesz fractional Complex Ginzburg-Landau-Schrödinger (CGLS) equation for various points of $x$ at $t = 0.5$ for $\alpha = 1.5$ and $\alpha = 1.75$.

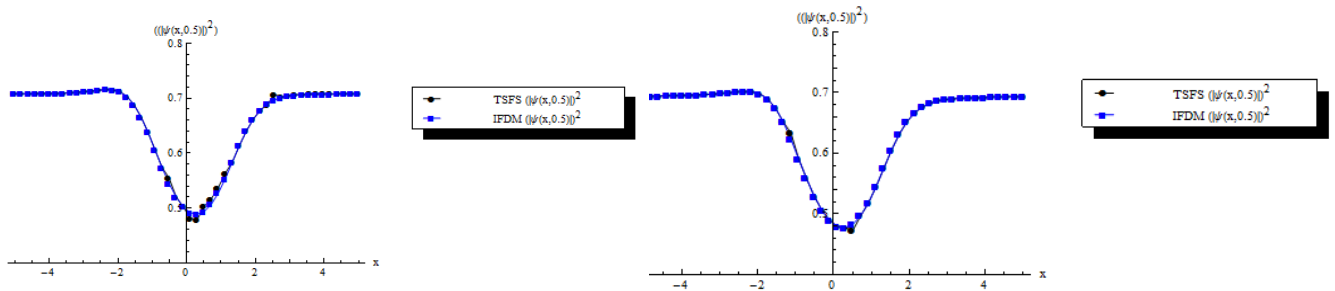

**Fig. 7.** Comparison of graphs for the solutions of $|\psi(x,t)|^2$ obtained from TSFS method and implicit finite difference(IFDM) scheme for the Riesz fractional Complex Ginzburg-Landau-Schrödinger (CGLS) equation (3.1) for Example 1 at $t = 0.5$ for $\alpha = 1.5$ and $\alpha = 1.75$.

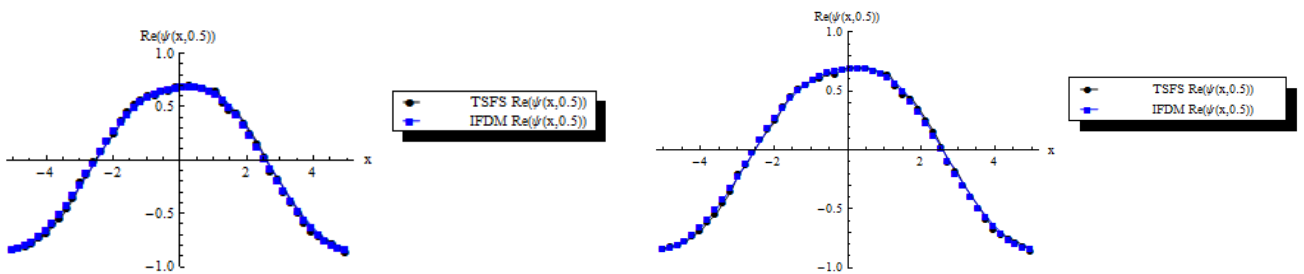

**Fig. 8.** Comparison of graphs for the solutions of Re($\psi(x,t)$) obtained from TSFS method and implicit finite difference(IFDM) scheme for the Riesz fractional Complex Ginzburg-Landau-Schrödinger (CGLS) equation (3.1) for Example 1 at $t = 0.5$ for $\alpha = 1.5$ and $\alpha = 1.75$.

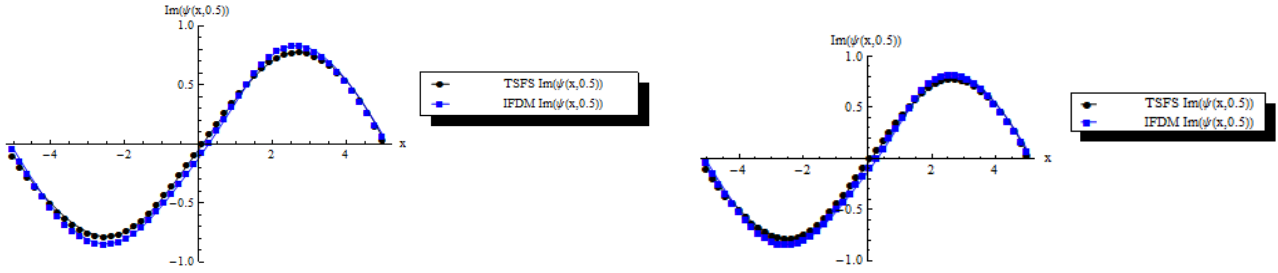

**Fig. 9.** Comparison of graphs for the solutions of $\text{Im}(\psi(x,t))$ obtained from TSFS method and implicit finite difference(IFDM) scheme for the Riesz fractional Complex Ginzburg-Landau-Schrödinger (CGLS) equation (3.1) for Example 1 at $t = 0.5$ for $\alpha = 1.5$ and $\alpha = 1.75$.

## 8. Conclusion

In the present paper the fractional CGLS equation with the Riesz fractional derivative has been tackled via an implicit second-order finite difference method (IFDM), where the fractional centered difference scheme is applied for the approximation of the Riesz fractional derivative. The stability of the implicit scheme has also been represented via von Neumann stability analysis. Further a comparative study has been put forth in the later section via a second order splitting scheme named Strang splitting technique for the purpose of investigating the efficiency of the results obtained. Moreover the spectral scheme is proved to be unconditionally stable. Further the tabulated results for absolute errors obtained by both the methods manifests that there is a fine tuning between the solutions obtained from both the approaches. The resulting values of the proposed methods have been utilized to display the error norms $L_2$ and $L_\infty$ for $|\psi(x,t)|^2$, $\text{Re}(\psi(x,t))$ and $\text{Im}(\psi(x,t))$ of Example 1. Furthermore, the results have also been displayed graphically.

### Acknowledgement

I am thankful to the authorities of my institute, National Institute of Technology to successfully carrying out this piece of research work and also all my professors and guidance from colleagues.